\def\timestamp{%
Time-stamp: <realcompact.tex: vrijdag 25-08-2023 at 17:02:36 (cest)>}
\def\stripname Time-stamp: <#1: #2 #3 at #4 #5>{#1: #3/#4}
\edef\filedate{\expandafter\stripname\timestamp}
\documentclass[a4paper]{amsart}

\DeclareMathSymbol\A0{AMSb}{`A}
\DeclareMathSymbol\D0{AMSb}{`D}
\DeclareMathSymbol\J0{AMSb}{`J}
\DeclareMathSymbol\K0{AMSb}{`K}
\DeclareMathSymbol\N0{AMSb}{`N}
\DeclareMathSymbol\R0{AMSb}{`R}
\DeclareMathSymbol\T0{AMSb}{`T}
\DeclareMathSymbol\V0{AMSb}{`V}
\DeclareMathSymbol\restr\mathbin{AMSa}{"16}

\newcommand\bee{\mathfrak{b}}
\newcommand\cee{\mathfrak{c}}
\newcommand\ess{\mathfrak{s}}
\newcommand\C{\mathfrak{C}}

\newcommand\calZ{\mathcal{Z}}
\newcommand\cl{\operatorname{cl}}
\newcommand\orpr[2]{\langle{#1},{#2}\rangle}

\theoremstyle{plain}
\newtheorem{proposition}{Proposition}[section]
\theoremstyle{remark}
\newtheorem{question}{Question}

\usepackage{amsrefs}

\begin{document}

\title{Some realcompact spaces}
\author[A. Dow]{Alan Dow}
\address{Department of Mathematics\\
         UNC-Charlotte\\
         9201 University City Blvd. \\
         Charlotte, NC 28223-0001}
\email{adow@charlotte.edu}
\urladdr{https://webpages.charlotte.edu/adow}

\author[K. P. Hart]{Klaas Pieter Hart}
\address{Faculty EEMCS\\TU Delft\\
         Postbus 5031\\2600~GA {} Delft\\the Netherlands}
\email{k.p.hart@tudelft.nl}
\urladdr{https://fa.ewi.tudelft.nl/\~{}hart}

\author[J. van Mill]{Jan van Mill}
\address{KdV Institute for Mathematics\\
         University of Amsterdam\\
         P.O. Box 94248\\
         1090~GE {} Amsterdam\\
         The Netherlands}
\email{j.vanmill@uva.nl}
\urladdr{https://staff.fnwi.uva.nl/j.vanmill/}

\author[J. Vermeer]{Hans Vermeer}
\address{Faculty EEMCS\\TU Delft\\
         Postbus 5031\\2600~GA {} Delft\\the Netherlands}
\email{j.vermeer@tudelft.nl}

\begin{abstract}
We present examples of realcompact spaces with closed subsets that 
are $C^*$-embedded but not $C$-embedded,
including one where the closed set is a copy of~$\N$.  
\end{abstract}

\date{\filedate}
\maketitle

\section*{Introduction}

The purpose of this note is to provide some examples of realcompact
(but not compact)
spaces that have closed subspaces that are $C^*$-embedded but not $C$-embedded,
and in particular an example where the closed subspace is a copy of the
discrete space~$\N$ of natural numbers --- 
what we henceforth call \emph{a closed copy of\/~$\N$}.

The reason for our interest is that we are not aware of any such examples.
The examples in~\cite{MR0407579}, for instance, of $C^*$- but not $C$-embedded
subsets are not all closed and when they are closed the pseudocompactness
of the ambient space makes $C$-embedding impossible.

The only explicit question of this nature that we could find
is in the paper~\cite{MR2367385}*{Question~1}, which asks whether 
$C^*$-embedded subsets (not necessarily closed) of first-countable spaces 
are $C$-embedded.
In that case there is an independence result: there is a counterexample
if $\bee=\ess=\cee$, and in the model obtained by adding supercompact
many random reals the implication holds, see~\cite{MR2823691}.

\smallskip  
The more specific question of having a closed copy of~$\N$ that is 
$C^*$-embedded but not $C$-embedded arises from an analysis of their
position in powers of the real line; see Section~\ref{sec:context}
for an explanation.

It is clear that our examples should be non-normal Tychonoff spaces.
After some preliminaries we briefly discuss two classical examples,
the Tychonoff and Dieudonn\'e planks, and introduce a further plank~$\J$.

The latter is pseudocompact but we modify it in two steps.
The first step yields a plank that is neither pseudocompact nor realcompact,
and the second step gives us our first example.

Our second example is constructed in Section~\ref{sec:another.plank}
and it contains a closed copy of~$\N$ that is $C^*$- but not $C$-embedded.

\section{Preliminaries}

We follow \cite{MR1039321} and \cite{MR0407579} as regards general topology
and rings of continuous functions.
As is common $C(X)$ and $C^*(X)$ denote the rings of real-valued continuous 
and bounded continuous functions respectively.

A subset~$A$ of a space~$X$ is \emph{$C$-embedded} if every continuous
function $f:A\to\R$ admits a continuous extension $\bar f:X\to\R$.
It is \emph{$C^*$-embedded} if every bounded continuous
function $f:A\to\R$ admits a bounded continuous extension $\bar f:X\to\R$.

We define a space $X$ to be \emph{realcompact} if it can be embedded into 
a power of the real line as a closed subset.
The most useful characterization for this paper is that every zero-set
ultrafilter with the countable intersection property has a non-empty
intersection, see~\cite{MR1039321}*{Theorem~3.11.11}.

\subsection*{Planks}

As noted above our examples will be non-normal Tychonoff spaces.
Non-normal because we need a closed subset that is not $C$-embedded
and Tychonoff because that is part of the definition of realcompactness.

There are various examples of such spaces, such as the Tychonoff plank~$\T$
(\cite{MR1512595} or~\cite{MR507446}*{Example~87}), and
the Dieudonn\'e plank~$\D$
(\cite{MR13297} or~\cite{MR507446}*{Example~89}).

Both start with the product set
$X=(\omega_1+1)\times(\omega_0+1)$ and take the subset
$P=X\setminus\{\orpr{\omega_1}{\omega_0}\}$ as the underlying set of the space.

In each case $P$~has the subspace topology where $X$~has a product
topology induced by topologies on the factors.

For~$\T$ one takes the order topologies on both ordinals.
For~$\D$ one enlarges the order topology of~$\omega_1+1$ by making all
points of~$\omega_1$ isolated.

We shall consider a third variation in Section~\ref{sec.J-variation} below.

\section{Context}\label{sec:context}

To begin we have the following proposition, which may be well-known 
but bears repeating here because it shows that if one has a 
non-$C$-embedded copy of~$\N$ in a realcompact space then that copy contains
many infinite subsets that \emph{are} $C$-embedded.

\begin{proposition}\label{prop.copy.of.N}
Let $X$ be realcompact and $A$ a subset whose closure is not compact,
then $A$~contains a countably infinite subset that is closed, discrete
and $C$-embedded in~$X$.  
\end{proposition}

\begin{proof}
Take a point~$x_0$ in $\beta X\setminus X$ that is in the closure of~$A$.
Apply~\cite{MR1039321}*{Theorem~3.11.10} to find a continuous function
$f:\beta X\to[0,1]$ such that $f(x_0)=0$ and $f(x)>0$ if $x\in X$.
Because $x_0$~is in the closure of~$A$ we can find a sequence 
$\langle a_n:n\in\N\rangle$ in~$A$ such that $\langle f(a_n):n\in\N\rangle$
is strictly decreasing with limit~$0$.

The set $N=\{a_n:n\in\N\}$ is closed and $C$-embedded in~$X$.
It is closed as a locally finite set of points.
If $g:N\to\R$ is given then we can take a continuous function $h:(0,1]\to\R$
such that $h(f(a_n))=g(a_n)$ for all~$n$. 
Then $h\circ f$ is a continuous extension of~$g$.
\end{proof}

The space in Section~\ref{sec:another.plank} illustrates this proposition
quite well: one can point out very many infinite $C$-embedded subsets
of the non-$C$-embedded copy of~$\N$ explicitly.  

This proposition also shows why the initial planks in 
Section~\ref{sec.J-variation} are not realcompact: there are not enough
$C$-embedded copies of~$\N$.

\subsection*{Closed copies of~$\N$ in other spaces}

Here we collect a few natural questions that arise when one considers
$C^*$- and $C$-embedding of closed copies of~$\N$.

Suppose one has two closed copies, $\N_1$ and $\N_2$ say, of the space
of natural numbers in a Tychonoff space~$X$.

\begin{enumerate}
\item If $\N_1$ and $\N_2$ are $C$-embedded is their union $C$-embedded?
      \label{questions}
\item If $\N_1$ and $\N_2$ are $C^*$-embedded is their union $C^*$-embedded?
\item If $\N_1$~is $C$-embedded and $\N_2$ is $C^*$-embedded is their 
      union $C^*$-embedded?
\end{enumerate}

Questions~(1) and~(3) have positive answers.
 
For question~(3) one uses a continuous extension $f:X\to\R$ of a 
bijection between~$\N_1$ and~$\N$ to obtain a discrete family $\{O_x:x\in\N_1\}$
of open sets with $x\in O_x$ for all~$x\in\N_1$.
Then, given a bounded function $g:\N_1\cup\N_2\to\R$ one first takes a bounded 
extension~$\bar g:X\to\R$ of~$g\restr\N_2$ and then modifies $\bar g$
on each~$O_x$ to obtain an extension of~$g$.

The argument for question~(1) is similar but easier because one can find
a single discrete family of open sets that separates the points of
$\N_1\cup\N_2$.

A counterexample to question~(2) can be obtained by taking Kat\v{e}tov's
example of a pseudocompact space with a closed $C^*$-embedded copy of~$\N$, 
see \cite{MR50264} or \cite{MR1039321}*{Example~3.10.29}.
The example is $\K=\beta\R\setminus\N^*$, and the copy of~$\N$ is 
just $\N$~itself. 
Take the sum of two copies of this space, $\K\times\{0,1\}$, and
for every $x\in\K\setminus\R$ identify the points $\orpr x0$ and~$\orpr x1$.
The copies $\N\times\{0\}$ and $\N\times\{1\}$ are both $C^*$-embedded
in the resulting quotient, but their union is not.

Below we shall show that question~(2) also has a negative answer in 
the class of realcompact spaces.

\subsection*{Closed copies of~$\N$ in powers of~$\R$}

The discrete space~$\N$ is realcompact, hence it admits many embeddings
into powers of~$\R$ as a closed and $C$-embedded set.

The specific question from the introduction is equivalent to the question
whether there is a closed copy of~$\N$ in some power of~$\R$ that 
is $C^*$-embedded but not $C$-embedded.
Indeed the latter is a special case of the former and a positive answer
to the former answers the latter by embedding the example as a closed
$C$-embedded copy into some power of~$\R$; the copy of~$\N$ is then
not $C$-embedded in that power.

The difference between $C$- and $C^*$-embedding manifests itself also
in the way certain maps can be factored through partial products.

\medbreak
Assume first that $\N$ is $C$-embedded in a power of~$\R$, say $\R^\kappa$.
Then there is a continuous function $f:\R^\kappa\to\R$ such that $f(n)=n$
for all $n\in\N$. 
It is well known that $f$~factors through a countable subset of~$\kappa$:
there are a countable subset~$I$ of~$\kappa$ and a continuous function
$g:\R^I\to\R$ such that $f=g\circ\pi$, where $\pi$ is the projection 
onto~$\R^I$, see~\cite{MR1039321}*{Problem~2.7.12}.
Then the projection~$\pi[\N]$ of~$\N$ in~$\R^I$ is $C$-embedded and we see
that every function from~$\N$ to~$\R$ has an extension that factors through
the partial power~$\R^I$.

\medbreak
Now assume $\N$ is $C^*$-embedded but not $C$-embedded in~$\R^\kappa$.
Then every bounded function from~$\N$ to~$[0,1]$ has a continuous extension
to~$\R^\kappa$.
Such a continuous extension will then factor through a partial product
with countably many factors but the set of factors will vary with the function.

\smallskip
Indeed, assume that there is a single countable set~$I$ such that every
bounded function $f:\N\to[0,1]$ has a continuous extension that factors
through~$\R^I$.
Apply this with the function $f(n)=2^{-n}$; using a factorization 
$\bar f=g\circ\pi$, as above, of an extension~$\bar f$ of~$f$ that 
the projection~$\pi$ onto~$\R^I$ is injective on~$\N$ and
that $\pi[\N]$~is relatively discrete in~$\R^I$.

We also find that $\pi[\N]$ is $C^*$-embedded in the metric space~$\R^J$,
and hence closed.
But then $\pi[\N]$ is $C$-embedded in~$\R^J$ and $\N$~is $C$-embedded 
in~$\R^\kappa$.
  
\bigskip
Using the plank~$\A$ from Section~\ref{sec:another.plank}
we obtain such a copy of~$\N$ in a power of~$\R$.
The standard embedding of~$\A$ in the power~$\R^{C(\A)}$ yields a closed
$C$-embedded copy of~$\A$.
The right-hand side~$R$ is a closed copy of~$\N$ that is $C^*$-embedded
in~$\A$ and hence in~$\R^{C(\A)}$, but not $C$-embedded in~$\R^{C(\A)}$.  

This then suggests the following question.

\begin{question}\label{question.1}
What is the minimum cardinal~$\kappa$ such that $\R^\kappa$ contains
a closed copy of~$\N$ that is $C^*$-embedded but not $C$-embedded?  
\end{question}

Since $\R^{\omega_0}$ is metrizable and, as we shall see, 
$\bigl|C(\A)\bigr|=\cee$ we know that $\aleph_0<\kappa\le\cee$.
This means that the Continuum Hypothesis settles this question, but
there may be some variation under other assumptions. 

Our answer to question~(2) from the list on page~\pageref{questions}
produces, in the same way, a closed copy of~$\N$ in~$\R^\cee$ that is 
not $C^*$-embedded.
After we submitted this paper we were able to answer the analogue
of Question~\ref{question.1}:
the smallest cardinal~$\kappa$ such that $\R^\kappa$ contains a closed 
copy of~$\N$ that is not $C^*$-embedded is~$\aleph_1$.
See~\cite{dhvmv2307.07223} for a surprising (to us) variety of closed copies 
of~$\N$ in~$\R^{\omega_1}$ that are not $C^*$-embedded.

\section{The plank $\J$ and a variation}\label{sec.J-variation}

In our third variation of the idea of the plank the topology on~$\omega_0+1$ 
remains as it is and we let, from now on, $\omega_1+1$~carry the topology of 
the one-point compactification of the discrete space~$\omega_1$, 
with $\omega_1$ the point at infinity.

In this case we denote the resulting space by~$\J$.
It is a minor variation of~\cite{MR1039321}*{Example~2.3.36}:
in the terminology of that book 
$\J=A(\aleph_1)\times A(\aleph_0)\setminus\{\orpr{x_0}{y_0}\}$, where
we have specified the underlying sets of the factors explicitly.

As in the case of $\T$ and $\D$ 
the top line $T=\omega_1\times\{\omega_0\}$ and 
the right-hand side $R=\{\omega_1\}\times\omega_0$ 
cannot be separated by open sets in~$\J$.
Hence their union is not $C^*$-embedded in the space~$\J$.

A more careful analysis of the continuous functions on~$\J$ will reveal that
neither~$T$ nor~$R$ is $C^*$-embedded.

Indeed: let $f:\J\to\R$ be continuous.
For each $n\in\omega_0$ the set 
$\{\alpha\in\omega_1:f(\alpha,n)\neq f(\omega_1,n)\}$
is countable.
It follows that there is an~$\alpha$ in~$\omega_1$ such that 
$f(\beta,n)=f(\omega_1,n)$ for all~$n$ and all  $\beta\ge\alpha$.
By continuity this implies that $f(\beta,\omega_0)=f(\alpha,\omega_0)$
for all $\beta\ge\alpha$. 
This shows that the function $\orpr\alpha{\omega_0}\mapsto\alpha\bmod2$
(the characteristic function of the odd ordinals),
which is continuous on~$T$, has no continuous extension to~$\J$.

If we let $r=f(\alpha,\omega_0)$ then it follows that 
$\lim_{n\to\infty}f(\omega_1,n)=r$.
We see that the function $\orpr{\omega_1}{n}\mapsto n\bmod 2$, which
is continuous on~$R$, has no continuous extension to~$\J$ either.

This argument also shows that $\J$ is not realcompact:
the co-countable sets on the top line form a zero-set ultrafilter
with the countable intersection property that has an empty intersection.
Alternatively use Proposition~\ref{prop.copy.of.N}: no infinite
subset of~$R$ is $C^*$-embedded.

The space~$\J$ is not pseudocompact either: the diagonal 
$\{\orpr nn:n\in\omega_0\}$ is a clopen discrete subset.

\subsection*{Ensuring $C^*$-embeddedness}

To ensure that $R$ is $C^*$-embedded we change the second factor in our product.

We let $X=(\omega_1+1)\times\beta\omega_0$ and
$P=X\setminus(\{\omega_1\}\times\omega_0^*)$.
The right-hand side~$R$ remains unchanged but the top line~$T$ now becomes
$\omega_1\times\omega_0^*$.

\smallbreak
To see why this makes the right-hand side $C^*$-embedded let 
$f:R\to[0,1]$ be continuous.
Take the unique continuous extension of $n\mapsto f(\omega_1,n)$ 
to~$\beta\omega_0$ and it on every vertical 
line~$\{\alpha\}\times\beta\omega_0$ to get an extension of~$f$ to the 
plank~$P$.

This does not make the right-hand side $C$-embedded: the analysis of the 
continuous functions on~$\J$ shows that for any extendable function~$f$
the function $n\mapsto f(\omega_1,n)$ should be extendable 
from~$\omega_0$ to~$\beta\omega_0$ and hence should be bounded.

When we adapt the analysis of continuous functions on~$\J$ to 
continuous functions on~$P$
we obtain that the intersection of a zero-set with the top line~$T$
contains a set of the form $A(\alpha,Z)=[\alpha,\omega_1)\times Z$, 
where $\alpha\in\omega_1$ and $Z$~is a zero-set of~$\omega_0^*$  
(and $Z$ could be empty of course).

Now take any point~$u$ in~$\omega_0^*$ and let $\calZ_u$ be the family
of zero-sets of $\omega_0^*$ that contain~$u$.
Then $\{A(\alpha,Z):\alpha\in\omega_1, Z\in\calZ_u\}$
generates a zero-set ultrafilter with the countable intersection property
that has an empty intersection.
Thus, the present plank is not realcompact.
Again, Proposition~\ref{prop.copy.of.N} applies as well: no closed copy
of~$\N$ (and there are many) in~$R$ is $C$-embedded.

\section{The plank $\V$}

It should be clear that the fact that continuous functions on~$\omega_1+1$
are constant on co-countable sets is the main cause that the two previous
examples are not realcompact.

To alleviate that we replace $\omega_1+1$ by $\beta\omega_1$, where $\omega_1$
still has the discrete topology.
We take the product $\Pi=\beta\omega_1\times\beta\omega_0$; our example
is $\V=\Pi\setminus (\omega_1^*\times\omega_0^*)$.

The top line and the right-hand side now become
$T=\omega_1\times\omega_0^*$ and
$R=\omega_1^*\times\omega_0$.

\subsection*{The right-hand side $R$ is $C^*$-embedded in~$\V$}

This is proved almost as in the case of the plank~$P$.

Let $f:R\to[0,1]$ be continuous.
Apply the Tietze-Urysohn extension theorem to each horizontal line~$H_n$
to obtain a continuous extension $f_n:H_n\to[0,1]$ of the restriction of~$f$
to $\omega_1^*\times\{n\}$.

Next take, for each $\alpha\in\omega_1$, the unique extension~$g_\alpha$
of the map $\orpr\alpha n\mapsto f_n(\alpha,n)$ 
to~$\{\alpha\}\times\beta\omega_0$.
The union of the maps $g_\alpha$ and~$f_n$ is an extension of~$f$ to~$\V$.

\subsection*{The right-hand side $R$ is not $C$-embedded in~$\V$}

Define $f:R\to\R$ by $f(x,n)=n$.
Assume $g:\V\to\R$ is a continuous extension of~$f$.
For each~$n$ and $k$ the set
$$
\{\alpha\in\omega_1:\bigl|g(\alpha,n)-n\bigr|\ge2^{-k}\}
$$
is finite, hence for each~$n$ the set $\{\alpha:g(\alpha,n)\neq n\}$
is countable.
It follows that there are co-countably many~$\alpha\in\omega_1$
such that $g(\alpha,n)=n$ for all~$n$.
For each such $\alpha$ the restriction of~$g$ to the compact 
set $\{\alpha\}\times\beta\omega_0$ would be unbounded, which is 
a contradiction.

\subsection*{The space $\V$ is realcompact}

Let $\calZ$ be a zero-set ultrafilter with the countable intersection
property.
We show that its intersection is nonempty.

To begin: if for some $n$ the clopen 
`horizontal line' $H_n=\beta\omega_1\times\{n\}$ belongs to~$\calZ$ then the
compactness of this line implies that $\bigcap\calZ$ is nonempty.

In the opposite case the complements of the $H_n$ belong to~$\calZ$;
the intersection of these complements is equal 
to the top line~$T$.
By the countable intersection property we find that every member of $\calZ$
intersects~$T$, hence $T\in\calZ$. 

For every subset $A$ of $\omega_1$ the partial top line 
$T_A=A\times\omega_0^*$ is a zero-set as it is the intersection of~$T$
with the clopen subset $\cl A\times\beta\omega_0$ of~$\Pi$.

Consider $u=\{A:T_A\in\calZ\}$.
This is an ultrafilter on~$\omega_1$ and it has
the countable intersection property and therefore, because $\omega_1$~is not
a measurable cardinal, it is a principal ultrafilter.
Take $\alpha\in\omega_1$ such that
$u=\{A\subseteq\omega_1:\alpha\in A\}$.

It follows that the compact set $\{\alpha\}\times\omega_0^*$ belongs to~$\calZ$
so that $\bigcap\calZ\neq\emptyset$.

\subsection*{Comments}

The natural maps from $\beta\omega_1$ onto~$\omega_1+1$
and and from~$\beta\omega_0$ 
onto~${\omega_0+1}$-as-the-one-point-compactification are perfect and 
irreducible.
Hence so is the product map from~$\Pi$ onto $(\omega_1+1)\times(\omega_0+1)$.
It follows that the restriction of this map to~$\V$ is perfect as well,
because $\V$~is the preimage of~$\J$.

We have seen that $\J$~is not realcompact, so we have here a very simple
perfect map that does not preserve realcompactness.

We also note that $\V$~is extremally disconnected and it is in fact
the absolute of~$\J$.

\section{Another plank}\label{sec:another.plank}

In this section we construct a realcompact space with a closed copy of~$\N$
that is $C^*$-embedded but not $C$-embedded.

We let $D$ be the tree $2^{<\omega}$ with the discrete topology and we 
topologize $D\cup2^\omega$ so as to obtain a natural compactification~$cD$ 
of~$D$.
If $x\in2^\omega$ then its $n$th neighbourhood~$U(x,n)$ will be the `wedge'
above $x\restr n$:
$$
U(x,n)=\{s\in cD: x\restr n\subseteq s\}
$$
Let $e:\beta D\to cD$ be extension of the identity map.

This yields a partition of~$D^*$ into closed sets, indexed by~$2^\omega$:
simply let $K_x=\{u\in D^*:e(u)=x\}$.

To construct our plank we take a point~$\infty$ not in~$2^\omega$ and 
topologize $\C=2^\omega\cup\{\infty\}$ by making every point of~$2^\omega$ 
isolated and letting 
$$
\{U:\infty\in U\land |\C\setminus U|\le\aleph_0\}
$$
be a local base at~$\infty$.

Let us note that $\C$~has a property in common with the horizontal lines
in our planks above: for every
continuous function $f:\C\to\R$ there is a neighbourhood of~$\infty$
(a co-countable set) on which $f$~is constant.

We let $\A$ be the following subspace of $\C\times\beta D$:
$$
\A=(\C\times D)\cup\bigcup_{x\in2^\omega}\{x\}\times K_x
$$
We let $R=\{\infty\}\times D$ denote the right-hand side of the plank.
The top line $T=\bigcup_{x\in2^\omega}\{x\}\times K_x$ is 
not as smooth as in the other examples; every point~$u$ of~$D^*$ occurs just
once in the top line, when $e(u)=x$.

\subsection*{$R$~is $C^*$-embedded}

This is as in the previous examples: $R$~is even $C^*$-embedded in
$R\cup(2^\omega\times\beta D)$.
Given $f:R\to[0,1]$ let $g:\beta D\to[0,1]$ be the \v{C}ech-Stone extension
of~$s\mapsto f(\infty,s)$ and then define $\bar f:\A\setminus R\to[0,1]$
by $\bar f(x,u)=g(u)$ (replicate $g$ on each vertical line but restrict it
to $\{x\}\times(\omega_0\cup K_x)$ each time).
Then $f\cup\bar f$ is a continuous extension of~$f$.

\subsection*{$R$~is not $C$-embedded} 

Below we show that $\A$~is realcompact, so Proposition~\ref{prop.copy.of.N}
implies that $R$~has many infinite $C$-embedded subsets.
Therefore the unbounded function without continuous extension must be chosen
with some care.

Define $f(\infty,s)=|s|$ (the length of~$s$).
Assume $g:\A\to\R$ is a continuous extension.
By the remark above there is a neighbourhood~$U$ of~$\infty$
such that $g$~is constant on $U\times\{s\}$ for every $s\in D$.
But then for every $x\in U\setminus\{\infty\}$ and $n\in\omega_0$
we have $g(x,x\restr n)=g(\infty,x\restr n)=f(\infty,x\restr n)=n$.
Since $K_x=\bigcap_n\cl_{\beta D}\{x\restr i:i\ge n\}$ this would imply
that $g(x,u)\ge n$ for all~$n$ when $u\in K_x$.

\subsection*{$\A$ is realcompact}

In the plank~$P$ in Section~\ref{sec.J-variation} we used~$\omega_0^*$
everywhere in the top line.
Combined with the fact that continuous functions were
constant on a tail on each horizontal line this implied that $P$~is
not realcompact, mainly because unbounded (to the right)
zero-sets in the top line contain sets of the form 
$[\alpha,\omega_1)\times Z$, where $Z$~is a zero-set of~$\omega_0^*$.
In the present example the disjointness of the~$K_x$ will provide us
with a richer supply of zero-sets; these will help ensure realcompactness
of~$\A$.  

Let $\calZ$ be a zero-set ultrafilter on~$\A$ with the countable
intersection property.

\smallbreak
For each $s\in D$ the horizontal $\C\times\{s\}$ is clopen, hence a zero-set.

The continuous function $f:\A\to[0,1]$ determined by setting $f(x,s)=2^{-|s|}$ 
for all $\orpr xs\in\C\times D$ has the top line~$T$ as its zero-set.

This means that we have a partition of~$\A$ into countably many zero-sets.
It follows that one of these sets must belong to~$\calZ$.

\smallbreak
If $\C\times\{s\}\in\calZ$ then either $\orpr\infty s\in\bigcap\calZ$
or there is a $Z\in\calZ$ is such that $\infty\notin Z$.
But then $Z$~is discrete and countable because
$\{x\in\C:\orpr xs\notin Z\}$ is open in~$\C$ and contains~$\infty$.
Then $\calZ$ determines a countably complete ultrafilter on~$Z$, which is fixed
because $|Z|$~is countable.

\smallbreak
We are left with the case that $T\in\calZ$.
Here is where we use the partition $\{K_x:x\in2^\omega\}$ of~$D^*$ to show that
$T$~may be split into zero-sets in many ways.

We show that whenever $A$~is clopen in the Cantor~set $2^\omega$ the 
union $Z(A)=\bigcup_{x\in A}\{x\}\times K_x$ is a zero-set in~$\A$.

By compactness and zero-dimensionality of~$cD$ we know there is a continuous 
function $f:cD\to\{0,1\}$ such that $f[A]=\{0\}$ 
and $f[2^\omega\setminus A]=\{1\}$ (we assume both $A$ and its complement are
non-empty).

We use $f$ to define $F:\A\to\{0,1\}$ by $F(x,s)=f(s)$  
if $\orpr xs\in\C\times D$ and $F(x,u)=f(x)$ if $u\in K_x$.

The function~$F$ is continuous on~$\A$ and we have
$Z(A)=T\cap Z_F$, so $Z(A)$~is a zero-set of~$\A$.

Using this we build countably many pairs of complementary zero sets
in~$T$.
For every $n\in\omega$ we let $A_n=\{x\in2^\omega:x(n)=0\}$
and $B_n=\{x\in2^\omega:x(n)=1\}$; these clopen sets determine
the zero-sets 
$Z(n,0)=\bigcup_{x\in A_n}\{x\}\times K_x$ and
$Z(n,1)=\bigcup_{x\in B_n}\{x\}\times K_x$ respectively.

Since $\calZ$ is a zero-set ultrafilter and $T\in\calZ$ we deduce that
for every~$n$ there is an element~$x(n)$ of~$\{0,1\}$ such that 
$Z(n,x(n))\in\calZ$.
Thus we get an~$x\in2^\omega$ such that $\{Z(n,x(n)):n\in\omega\}$ is a 
subfamily of~$\calZ$.

Its intersection is equal to~$\{x\}\times K_x$ and because $\calZ$
has the countable intersection property this compact set belongs to~$\calZ$.
It follows that $\bigcap\calZ\neq\emptyset$.

\medskip
As mentioned before, Proposition~\ref{prop.copy.of.N} implies that $R$~has many
infinite $C$-em\-bed\-ded subsets.
A lot of these can be pointed out explicitly.

For every $x\in2^\omega$ the set $N_x=\{\orpr\infty{x\restr n}:n\in\omega\}$
is $C$-embedded in~$\A$.
Given a function $f:N_x\to\R$ we extend it to~$R$ first by setting 
$\bar f(\infty,s)=0$ for all other~$s$.
Then we extend $\bar f$ horizontally: $\bar f(y,s)=\bar f(\infty,s)$ for 
all~$y$ and~$s$, \emph{except} for $y=x$, we set $\bar f(x,s)=0$ for all~$s$.
Now we can set $\bar f(t)=0$ for all $t$ in the top line to get our continuous
extension to all of~$\A$.

In a similar fashion every infinite antichain in~$2^{<\omega}$ yields 
an infinite $C$-embedded subset as well.

\subsection*{More answers}

We can use $\A$ and some variations to answer some of the questions raised
earlier in this paper. 

\subsubsection*{The smallest power of~$\R$}
The set $\C\times D$ is dense in~$\A$, so every member of~$C(\A)$ is determined
by its restriction to this set.
Using the fact that continuous functions on~$\C$ are constant on co-countable 
sets we see that there are $\cee$~many such restrictions.
We conclude that $C(\A)$~has cardinality~$\cee$, as claimed in the discussion
of Question~\ref{question.1}.

\subsubsection*{The union of two closed $C^*$-embedded copies of~$\N$}
We can use $\A$ much like we used~$\K$ to create a realcompact space
with two closed $C^*$-embedded copies of~$\N$ whose union is not $C^*$-embedded.
Take $\A\times\{0,1\}$ and identify the points~$\orpr t0$ and~$\orpr t1$ for
all~$t$ in the top line~$T$.
Then $R\times\{0\}$ and $R\times\{1\}$ are still $C^*$-embedded in the 
resulting quotient space, but their union is not: mapping $\orpr ri$ to~$i$
results in a bounded function without a continuous extension.
The proof that the quotient space is realcompact is almost verbatim that
of the realcompactness of~$\A$. 
Note that the $R\times\{0\}$ and $R\times\{1\}$ are separated 
(neither intersects the closure of the other), so their union
is a closed copy of~$\N$ that is not $C^*$-embedded.
The quotient space also has $\cee$~many real-valued continuous functions,
hence also we obtain a closed copy of~$\N$ in~$\R^\cee$ that is not 
$C^*$-embedded.
This copy is quite unlike the closed copies of~$\N$ in~$\R^{\omega_1}$
that are constructed in~\cite{dhvmv2307.07223}.

\subsubsection*{Another closed copy of~$\N$ that is not $C^*$-embedded}
If we replace $\beta D$ by~$cD$ in~$\A$ then we obtain a realcompact plank
where the right-hand side is a closed copy of~$\N$ that is not $C^*$-embedded.

The analogue of~$\A$ is the following subspace of $\C\times cD$:
$$
(\C\times D)\cup\{\orpr xx:x\in2^\omega\}
$$
That this space is realcompact is shown exactly as for~$\A$.
However in this space the right-hand side~$R$ is not $C^*$-embedded.

Since $2^\omega$ is homeomorphic to its own square it is relatively easy
to produce two disjoint open sets $U$ and $V$ in~$2^\omega$ with a dense union
and whose common boundary~$F$ is homeomorphic to~$2^\omega$ itself.

Via the map $e:\beta D\to cD$ we can find a subset $C$ of $D$ such that
$\cl U\subseteq\cl C$ and $\cl V\subseteq\cl(D\setminus C)$.

Define $f:R\to[0,1]$ by $f(\infty,s)=\chi(s)$, where $\chi$ is the 
characteristic function of~$C$.
As before, given a continuous extension~$\bar f$ of~$f$, we would have a 
countable set~$B$ such that $\bar f(x,s)=f(\infty,s)$ for 
all~$x\in2^\omega\setminus B$ and all $s\in D$.
But then $\bar f$ would not be continuous at~$\orpr xx$ 
whenever $x\in F\setminus B$.

\end{document}